\documentclass[10pt]{article}
\usepackage{amsfonts,color}
\usepackage{amsmath,amssymb}
\usepackage{epsfig}
\usepackage{lscape}
\usepackage{amssymb}
\usepackage{graphicx}
\usepackage[utf8]{inputenc}	
\usepackage{amsthm,amssymb}
\usepackage{amsfonts}
\usepackage[english]{babel}
\usepackage{dsfont}
\usepackage{bbm}
\usepackage[T1]{fontenc}
\usepackage{colonequals}
\usepackage{geometry}
\usepackage{mathtools}
\usepackage{csquotes}
\usepackage{url}
\usepackage{tikz}
\usetikzlibrary{arrows}

 \newtheorem{theorem}{Theorem}[section]

 \newtheorem{definition}[theorem]{Definition}

\DeclareMathOperator{\tr}{tr}

\newcommand{\innerproduct}[1]{\langle #1 \rangle}
\newcommand{\iprod}{\innerproduct}
\newcommand{\norm}[1]{\lVert #1 \rVert}
\newcommand{\Cof}{{\rm Cof}}

\makeatletter
\let\@fnsymbol\@arabic
\makeatother

\newcommand{\id}{{\boldsymbol{\mathbbm{1}}}}

\newcommand{\setvert}{\:|\:}

\newcommand{\edit}[1]{#1} 
\newcommand{\rev}[1]{#1} 
\newcommand{\revtwo}[1]{#1} 

\begin{document}
\title{Injectivity of the Cauchy-stress tensor along rank-one connected lines under  strict rank-one convexity condition}

\author{
Patrizio Neff\thanks{	Corresponding author: Patrizio Neff,   \ Head of Lehrstuhl f\"{u}r Nichtlineare Analysis und Modellierung, Fakult\"{a}t f\"{u}r
		Mathematik, Universit\"{a}t Duisburg-Essen,  Thea-Leymann Str. 9, 45127 Essen, Germany, email: patrizio.neff@uni-due.de} \quad and\quad L. Angela Mihai\thanks{L. Angela Mihai, Senior Lecturer in Applied Mathematics, School of Mathematics, Cardiff University, Seghennydd Road, Cardiff, CF24 4AG, UK, \ email: MihaiLA@cardiff.ac.uk}
}
\date{\today\vspace*{-1em}}
\newgeometry{top=2.1em}
\maketitle
\begin{abstract}
	In this note, we show that the Cauchy stress tensor $\sigma$ in nonlinear elasticity is injective along rank-one connected lines provided that the constitutive law is strictly rank-one convex. This means that $\sigma(F+\xi\otimes \eta)=\sigma(F)$ implies $\xi\otimes \eta=0$ under strict rank-one convexity.
	\rev{As a consequence of this seemingly unnoticed observation, it follows that rank-one convexity and a homogeneous Cauchy stress \revtwo{imply} that the left Cauchy-Green strain is homogeneous, as is shown in \cite{Mihai_Neff_homogeneous}.}
\\[1.4em]
\\[1.4em]
\end{abstract}
\tableofcontents
%
%
\newgeometry{top=8.4em}

\section{Introduction}

The search for a priori constitutive inequalities has been termed by Truesdell \cite{Truesdell56,Truesdell60} the ``Hauptproblem" of nonlinear elasticity. These constitutive inequalities  should  guarantee reasonable physical response under all possible circumstances \cite[18.6.3]{Silhavy_book}. We focus here on one of these requirements, namely rank-one convexity, and exhibit a hitherto \rev{unnoticed} consequence of strict rank-one convexity in connection with the Cauchy stress tensor. \\

Following a definition by Ball \cite[Definition 3.2]{Ball77}, we say that $W$ is  \emph{strictly rank-one convex} on ${\rm GL}^+(3)=\{X\in\mathbb{R}^{3\times 3} \,\setvert\,$ $ \det X > 0\}$ if it is strictly convex on all closed line segments in ${\rm GL}^+(3)$ with end points differing by a matrix of rank one, \rev{i.e.,}
\begin{align}
W( F+(1-\theta)\, \xi\otimes \eta)< \theta \,W( F)+(1-\theta) W(F+\xi\otimes \eta)
\end{align}
for all $F\in {\rm GL}^+(3)$, $\theta\in[0,1]$ and all $\,\, \xi,\, \eta\in\mathbb{R}^3$ with $F+t\, \xi\otimes \eta\in {\rm GL}^+(3)$ for all
$t\in[0,1]$, where $\xi\otimes\eta$ denotes the dyadic product. \rev{Rank-one} convexity is connected to the study of wave propagation 
\cite{SawyersRivlin78} or hyperbolicity of the dynamical equations of elasticity, and plays an important role in the existence and uniqueness theory for linear elastostatics and elastodynamics \cite{Ogden83,fosdick2007note,edelstein1968note,simpson2008bifurcation}, cf.\ \cite{knowles1976failure}.
 Important criteria for the rank-one convexity of \rev{stored energy density} functions were first established by Knowles and Sternberg \cite{knowles1975ellipticity}, see also \cite{MartinGhibaNeff,NeffGhibaLankeit,GhibaMartinNeff}.

 In this paper we use  the Frobenius tensor norm
 $\|{X}\|^2=\langle {X},{X}\rangle_{\mathbb{R}^{n\times n}}$, where $\langle {X},{Y}\rangle_{\mathbb{R}^{n\times n}}$ is the standard Euclidean scalar product on $\mathbb{R}^{n\times n}$. If no  confusion can arise, we will suppress the subscripts $\mathbb{R}^{n\times n}$. The identity tensor on $\mathbb{R}^{n\times n}$ will be denoted by $\id$, so that
 $\tr{(X)}=\langle {X},{\id}\rangle$.\\
 
Rank-one convexity is preferably expressed in terms of the \rev{stored energy} density $W(F)$ or as a monotonicity requirement for the first Piola-Kirchhoff stress tensor $S_1=D\, W(F)$ along rank-one lines, \revtwo{i.e.,}
\begin{align}\label{sroc}
\langle S_1(F+\xi\otimes \eta)-S_1(F), \xi \otimes \eta\rangle_{\mathbb{R}^{3\times 3}} >0\qquad \forall\ \xi\otimes \eta\neq 0, \qquad \forall \ F\in {\rm GL}^+(3),
\end{align}
which, if $W$ is twice-differentiable, turns into the well-known strong-ellipticity condition
\begin{align}
D^2W(F)(\xi\otimes \eta, \xi\otimes \eta)>0 \qquad \forall\,\xi\otimes \eta\neq 0, \qquad \forall \ F\in {\rm GL}^+(3).
\end{align}

Since objective \rev{stored energy density functions} cannot be convex in $F$ \cite{Schroeder_Neff_IJSS}, the first Piola-Kirchhoff stress $S_1(F)$ will, in general, not be injective (\cite[6.2.38]{Ogden83},\cite[18.4.5]{Silhavy_book}). However, the strict monotonicity condition \eqref{sroc} means that
$
 S_1(F+\xi\otimes \eta)=S_1(F)
$ implies $\xi\otimes \eta=0$. This motivates the following

\begin{definition} The stress tensor $S$ is rank-one injective at $F$ if
	\begin{align}\label{eg23}
	S(F+\xi\otimes \eta)=S(F)\qquad\Longleftrightarrow\qquad\xi\otimes \eta= 0\,.
	\end{align}
\end{definition}

In this sense, if the \rev{stored energy density} is strictly rank-one convex, then the first Piola-Kirchhoff stress tensor $S_1(F)$ is everywhere rank-one injective.\\

The only well-known consequence of rank-one convexity in connection to the Cauchy stress tensor are the Baker-Ericksen inequalities \cite{Baker_Ericksen} for the principal \rev{values of the} Cauchy stress. These, however, are meaningful only for isotropy \cite{Fosdick_Silhavy}.
 
Here, we show by a short and elementary calculation that strict monotonicity of the first Piola-Kirchhoff stress tensor $S_1$ along rank-one lines implies injectivity of the Cauchy stress tensor along rank-one lines.

\rev{This elementary observation} answers a question raised \edit{in a recent contribution} \cite{Mihai_Neff_homogeneous}: \rev{Is it} impossible for a strictly rank-one convex \rev{stored energy to admit a continuous deformation that corresponds to a homogeneous Cauchy stress field but has jumps in its deformation gradient field across planar interfaces?}
Indeed, in \cite{Mihai_Neff_homogeneous} we \rev{show} that a non rank-one convex formulation may allow for \rev{a deformation with a} homogeneous Cauchy stress field but \rev{an} inhomogeneous \rev{left Cauchy-Green strain field}.\\

\rev{%
	We consider the following general situation: Let
	\[
		\sigma\colon{\rm GL}^+(3)\to{\rm Sym}(3)\,,\quad F\mapsto\sigma(F)
	\]
	denote the Cauchy stress response function induced by the stored energy density $W$, and let $F\in{\rm GL}^+(3)$ \revtwo{be such that}
	\begin{equation}\label{eg2}
		\sigma(F+\xi\otimes \eta)=\sigma(F)
	\end{equation}
	for some $\xi\otimes 
	\eta\neq 0$. We recall the basic relation \cite{Ciarlet}
	\begin{equation}\label{cKr}
		\sigma(F)=S_1(F)\, ({\rm Cof}(F))^{-1}
	\end{equation}
	and note that in case of isotropy we may write 
	\begin{align}
		\sigma(F)&=\widetilde{\sigma}(F\, F^T)=\widetilde{\sigma}(B)\, , \quad\quad 
		\widetilde{\sigma}\colon{\rm Sym}^+(3)\to {\rm Sym}(3),\quad B\mapsto \widetilde{\sigma}(B)\,.
	\end{align}
}

In isotropic nonlinear elasticity, a number of energies (suitable Neo-Hooke, Mooney-Rivlin \cite{Ciarlet,Ogden83}, the exponentiated Hencky energy \cite{NeffGhibaLankeit}) define an invertible Cauchy stress-strain relation, in the sense that the mapping  $B\mapsto \widetilde{\sigma}(B)$ is invertible. In this case $\sigma(F+\xi\otimes \eta)=\widetilde{\sigma}(\widehat{B})=\widetilde{\sigma}(B)=\sigma(F)$ leads to $B=\widehat{B}$. This, together with ${\rm det}\,\widehat{F}={\rm det}\,F>0$ implies $\xi\otimes \eta=0$ in \eqref{eg2}. A \edit{self-contained elementary} proof of this fact is given in the appendix.

Our subsequent development  will be independent of any invertibility assumption for the Cauchy stress $\sigma$ in the isotropic representation with $\widetilde{\sigma}$. 

\section{Injectivity of the Cauchy-stress tensor along rank-one lines for strictly rank-one convex energies}\setcounter{equation}{0}

We will show that \edit{equality} \eqref{eg2} combined with  strict rank-one convexity in the format of \eqref{sroc} leads to a  contradiction.\footnote{\edit{%
	The following alternative proof, which uses the identity ${\rm Cof}(F+\xi\otimes\eta).\eta = {\rm Cof}F.\eta$, see \cite[eq.~1.1.18]{Silhavy_book}, was \edit{kindly} suggested \revtwo{by the reviewer}:
	\begin{align*}
		\iprod{S_1(F+\xi\otimes\eta) - S_1(F),\, \xi\otimes\eta} &= \iprod{\sigma(F+\xi\otimes\eta) \cdot \Cof(F+\xi\otimes\eta) - \sigma(F)\cdot \Cof F,\, \xi\otimes\eta}\\
		&= \iprod{\sigma(F+\xi\otimes\eta),\, (\xi\otimes\eta)\big(\Cof(F+\xi\otimes\eta)\big)^T} - \iprod{\sigma(F),\, (\xi\otimes\eta)(\Cof F)^T}\\
		&= \iprod{\sigma(F+\xi\otimes\eta),\; \xi\otimes(\Cof(F+\xi\otimes\eta).\eta)} - \iprod{\sigma(F),\, \xi\otimes(\Cof F.\eta)}\\
		&= \iprod{\sigma(F+\xi\otimes\eta),\; \xi\otimes(\Cof(F).\eta)} - \iprod{\sigma(F),\, \xi\otimes(\Cof F.\eta)}\\
		&= \iprod{\sigma(F+\xi\otimes\eta) - \sigma(F),\, \xi\otimes((\Cof F).\eta)}\,.
	\end{align*}
	If the \revtwo{stored energy density function} is strictly rank-one convex, the latter identity implies that if $\sigma(F+\xi\otimes\eta)=\sigma(F)$, then $\xi\otimes\eta=0$.
}}
\begin{proof}
To this aim, using \eqref{cKr} we compute
\begin{align}
\sigma(F+\xi\otimes \eta)=S_1(F+\xi\otimes \eta)\,({\rm Cof}(F+\xi\otimes \eta))^{-1},\qquad \sigma(F)=S_1(F)\, ({\rm Cof}(F))^{-1}.
\end{align}
Hence, from \eqref{eg2} it follows that
\begin{alignat}{2}
&&S_1(F+\xi\otimes \eta)({\rm Cof}(F+\xi\otimes \eta))^{-1}&=S_1(F)\, ({\rm Cof}(F))^{-1}\notag\\
\iff \quad&&S_1(F+\xi\otimes \eta)&=S_1(F)\, ({\rm Cof}(F))^{-1}\,({\rm Cof}(F+\xi\otimes \eta))\,. \label{eq:PKoneFplusRankOne}
\end{alignat}
Since ${\rm Cof }(A)\, {\rm Cof }(B)={\rm Cof} (A\, B)$ and $ ({\rm Cof}A)^{-1}={\rm Cof}(A^{-1})$ for all $A, B\in {\rm GL}^+(3)$, we obtain
\begin{align}
S_1(F+\xi\otimes \eta)=S_1(F)\, {\rm Cof}\,(F^{-1}F+F^{-1}\xi\otimes \eta)=S_1(F)\, {\rm Cof}\,(\id+F^{-1}\xi\otimes \eta)\,.
\end{align}
Using now the expansion
${\rm Cof}\,(\id+H)={\rm Cof}\,(\id)+D\,{\rm Cof}(F)\Big|_{\id}.\,H+{\rm Cof}\,(H)$, see \cite{Schroeder_Neff_poly}, we find
\begin{align}
{\rm Cof}\,(\id+F^{-1}\xi\otimes \eta)={\rm Cof}\,(\id)+D\,{\rm Cof}(F)\Big|_{\id}.\,(F^{-1}\xi\otimes \eta)+\underbrace{{\rm Cof}\,(F^{-1}\xi\otimes \eta)}_{=0}\,,
\end{align}
and since
\begin{align}
D\,{\rm Cof}(F).\, H=\left(\langle F^{-T}, H\rangle \, \id-F^{-T} H^T\right){\rm Cof} F\ \Rightarrow\  D\,{\rm Cof}(F)\Big|_{\id}.\, H=\langle \id, H\rangle \, \id- H^T\,,
\end{align}
we can rewrite equality \eqref{eq:PKoneFplusRankOne} as
\begin{align}
S_1(F+\xi\otimes \eta)&=S_1(F)\, \left[ \id+ D\,{\rm Cof}(F)\Big|_{\id}.\,(F^{-1}\xi\otimes \eta)\right]\\
&=S_1(F)\, \left[ \id+\langle \id, (F^{-1}\xi\otimes \eta)\rangle \, \id- (F^{-1}\xi\otimes \eta)^T\right].\notag
\end{align}
Going back to the strict rank-one convexity condition \eqref{sroc}, we compute now
\begin{align}
\langle S_1(F+\xi\otimes \eta)-S_1(F), \xi\otimes \eta\rangle &= \langle S_1(F)\, \left[ \id+\langle \id, (F^{-1}\xi\otimes \eta)\rangle \, \id- (F^{-1}\xi\otimes \eta)^T\right]-S_1(F), \xi\otimes \eta \rangle\notag \\
&=\langle  \langle \id, (F^{-1}\xi\otimes \eta)\rangle \, S_1(F)\,- S_1(F)\,(F^{-1}\xi\otimes \eta)^T, \xi\otimes \eta \rangle
\notag \\
&=\langle \id, (F^{-1}\xi\otimes \eta)\rangle \,\langle   S_1(F), \xi\otimes \eta \rangle\,-\langle  S_1(F),\, (\xi\otimes \eta) \, \,(F^{-1}\xi\otimes \eta)\rangle
\notag \\
&=\langle \id, (F^{-1}\xi\otimes \eta)\rangle \,\langle   S_1(F), \xi\otimes \eta \rangle\,-\langle  S_1(F),\, \langle\eta ,F^{-1}\xi\rangle\, (\xi \otimes \eta)\rangle
\notag \\
&=\langle\eta ,F^{-1}\xi\rangle\,\langle   S_1(F), \xi\otimes \eta \rangle\,-\langle\eta, F^{-1}\xi\rangle\,\langle  S_1(F),  \xi \otimes \eta)\rangle=0\,.
\end{align}
Here, we have used that $\langle \id, a\otimes b\rangle_{\mathbb{R}^{3\times3}}=\langle b,a\rangle_{\mathbb{R}^3}$ as well as $(a\otimes b)\, (c\otimes d)=\langle b,c\rangle\,(a\otimes d)$, for all $a,b,c,d\in \mathbb{R}^3$.

Therefore, the assumption of the non-injectivity along rank-one lines \eqref{eg2} is in contradiction to the strict rank-one convexity \eqref{sroc}.
\end{proof}

\edit{In summary, we} have shown that  strict rank-one convexity implies that
\begin{align}\label{eg22}
\sigma(F+\xi\otimes \eta)=\sigma(F)\ \ \ \  \text{is impossible for a nontrivial}\ \ \ \  \xi\otimes 
\eta\neq 0, \ \ \ \xi,\eta\in \mathbb{R}^3. 
\end{align}

In these terms, we have thus proved that
\begin{align}
\text{strict rank-one convexity}\qquad \Longrightarrow\qquad &\text{the Cauchy stress} \quad \sigma\\ &\text{is rank-one injective for all}\ \ \ F\in{\rm GL}^+(3).\notag
\end{align}

\section{Conclusion}
Our simple calculation shows that for strictly rank-one convex \rev{stored energy density functions} it is impossible to have a constant Cauchy-stress field in response to a rank-one connected laminate microstructure. Our result suggests also that some form of injectivity for the Cauchy stress is natural to require in nonlinear elasticity and this injectivity should be the object of further studies.

\edit{In order to \rev{give added perspective to} our result on injectivity of the Cauchy stress, let us consider the uni-constant Blatz-Ko \rev{stored energy density function}
\[
	W(F) = \frac\mu2 \,\Big(\norm{F}^2 + \frac{2}{\det F} - 5\Big)\,.
\]
This function is strictly polyconvex, hence strictly rank-one elliptic with Cauchy stress
\begin{equation}
\label{eq:BKCauchy}
	\widetilde{\sigma}\colon {\rm Sym}^+(3) \to {\rm Sym}(3)\,,\qquad \widetilde{\sigma}(B) = \frac{\mu}{\det B}\,\Big(\sqrt{\det B} \cdot B - \id\Big)\,.
\end{equation}
The Cauchy stress in \eqref{eq:BKCauchy} is not bijective, which can be seen along the family $B=\alpha\cdot\id$, $\alpha>0$. The spherical part $\frac13\tr(\sigma)$ of the Cauchy stress first increases for increasing $\alpha$ and then decreases. Thus strict polyconvexity alone is not enough to prevent this unphysical response \cite{jog2013conditions}. We need to require a condition beyond polyconvexity. Injectivity of the Cauchy stress is a candidate implying the classical pressure-compression inequality \cite{NeffGhibaLankeit}
\begin{equation}
\frac 13\,\tr(\sigma(\lambda\,\id))\cdot[\lambda-1] > 0\,,
\end{equation}
which would already exclude the deficiency of the Blatz-Ko strain energy.}

\section{Acknowledgements}
The support for L.~Angela Mihai  by the Engineering and Physical Sciences Research Council of Great Britain under research grant EP/M011992/1 is gratefully acknowledged. \edit{We thank the reviewer for pointing out \rev{the shorter proofs noted in Section 2 and in the Appendix}.}

\bibliographystyle{plain} 

\begin{thebibliography}{10}
	

\bibitem{Baker_Ericksen}
	M. Baker and J.E. Ericksen
	\newblock Inequalities restricting the form of the stress-deformation relations for isotropic elastic solids and Reiner-Rivlin fluids.
		\newblock {\em  J. Wash. Acad. Sci. }, 44:33--35, 1954.


	\bibitem{Ball77}
	J.~M. Ball.
	\newblock Convexity conditions and existence theorems in nonlinear elasticity.
	\newblock {\em Arch. Rat. Mech. Anal.}, 63:337--403, 1977.
	
	\bibitem{Ciarlet} P.G. Ciarlet.
\newblock {\em Three-{D}imensional {E}lasticity}, Elsevier, 
	Studies in {M}athematics and its {A}pplications, Amsterdam, 
	1988.
	
	\bibitem{edelstein1968note}
	W.~Edelstein and R.~Fosdick.
	\newblock A note on non-uniqueness in linear elasticity theory.
	\newblock {\em Z. Angew. Math. Phys.}, 19(6):906--912, 1968.
	
	\bibitem{Fosdick_Silhavy}
	R.~Fosdick and M. Silhavy.
	\newblock Generalized Baker-Ericksen inequalities.
	\newblock {\em J. Elasticity}, 85:39--44, 2006.


	\bibitem{fosdick2007note}
	R.~Fosdick, M.D. Piccioni, and G.~Puglisi.
	\newblock A note on uniqueness in linear elastostatics.
	\newblock {\em J. Elasticity}, 88(1):79--86, 2007.
	
	\bibitem{GhibaMartinNeff}
	I.D. Ghiba,	R.J. Martin, and P.~Neff.
	\newblock Rank-one convexity implies polyconvexity in planar objective, isotropic and incompressible nonlinear elasticity.
	\newblock {\em submitted}, 2016.
\edit{		
	\bibitem{jog2013conditions}
	C.S.~Jog and K.D.~Patil.
	\newblock Conditions for the onset of elastic and material instabilities in hyperelastic materials.
	\newblock {\em Archive of Applied Mechanics}, 83.5:661--684, 2013.
}
	\bibitem{knowles1975ellipticity}
	J.K. Knowles and E.~Sternberg.
	\newblock On the ellipticity of the equations of nonlinear elastostatics for a
	special material.
	\newblock {\em J. Elasticity}, 5(3-4):341--361, 1975.
	
	\bibitem{knowles1976failure}
	J.K. Knowles and E.~Sternberg.
	\newblock On the failure of ellipticity of the equations for finite
	elastostatic plane strain.
	\newblock {\em Arch. Rat. Mech. Anal.}, 63(4):321--336, 1976.
	
	\bibitem{MartinGhibaNeff}
	R.J. Martin, I.D. Ghiba, and P.~Neff.
	\newblock Rank-one convexity implies polyconvexity for isotropic, objective and
	isochoric elastic energies in the two-dimensional case.
	\newblock {\em to appear in Proc. Roy. Soc. Edinburgh Sect. A}, 2016.
	

\bibitem{Mihai_Neff_homogeneous}
	L. A. Mihai and P.~Neff.
	\newblock Hyperelastic bodies under homogeneous Cauchy stress induced by non-homogeneous finite deformations.
	\newblock {\em to appear in Int. J. Nonl. Mechanics}, 2016.

	
	\bibitem{NeffGhibaLankeit}
	P.~Neff, I.~D. Ghiba, and J.~Lankeit.
	\newblock The exponentiated {H}encky-logarithmic strain energy. {P}art {I}:
	{C}onstitutive issues and rank--one convexity.
	\newblock {\em J. Elasticity}, 121:143--234, 2015.
	
	\bibitem{Ogden83}
	R.W. Ogden.
	\newblock {\em Non-{L}inear {E}lastic {D}eformations.}
	\newblock Mathematics and its Applications. Ellis Horwood, Chichester, 1.
	edition, 1983.
	
	\bibitem{SawyersRivlin78}
	K.N. {Sawyers} and R.~{Rivlin}.
	\newblock {On the speed of propagation of waves in a deformed compressible.
		elastic material.}
	\newblock {\em {Z. Angew. Math. Phys.}}, 29:245--251, 1978.
	
	
	
	\bibitem{simpson2008bifurcation}
	H.~Simpson and S.~Spector.
	\newblock On bifurcation in finite elasticity: buckling of a rectangular rod.
	\newblock {\em J. Elasticity}, 92(3):277--326, 2008.

\bibitem{Silhavy_book}
	M. Silhavy.
	\newblock The Mechanics and Thermodynamics of Continuous Media.
	\newblock {\em Springer}, 1997.	


\bibitem{Schroeder_Neff_IJSS}
	J. Schr\"oder and P. Neff.
	\newblock Poly-, Quasi- and Rank-One Convexity in Applied Mechanics.
	\newblock {\em Springer}, 2010.	

\bibitem{Schroeder_Neff_poly}
	J. Schr\"oder and P. Neff.
	\newblock Invariant formulation of hyperelastic transverse isotropy based on polyconvex free energy functions.
		\newblock {\em Int. J. Solids Struct.}, 40(2):401-445, 2003.

	\bibitem{Truesdell56}
		C.~Truesdell.
		\newblock Das ungel\"oste Hauptproblem der endlichen Elastizitätstheorie. \newblock {\em Z.  Angew. Math. Mech.}, 36: 97--103, 1956.
		
	\bibitem{Truesdell60}
	C.~Truesdell and R.~Toupin.
	\newblock The Classical Field Theories.
	\newblock In S.~Fl\"ugge, editor, {\em Handbuch der {P}hysik}, volume III/1.
	Springer, Heidelberg, 1960.
		
	
\end{thebibliography}
\addcontentsline{toc}{section}{References}

\begin{footnotesize}

\section{Appendix}

In this appendix we show\footnote{\edit{%
	\rev{The following alternative proof was kindly suggested by the reviewer: Rewriting} \eqref{eq_six} as
	\[
		(F\eta + \norm{\eta}^2\xi)\otimes\xi = - \xi\otimes F\eta
	\]
	and recalling that $a\otimes b = c\otimes d \neq 0$ if and only if there is $\lambda\in\mathbb{R}\setminus\{0\}$ such that $a=\lambda\,c$ and $b=\frac{d}{\lambda}$, then the assumption $\xi\otimes\eta\neq0$ implies
	\[
		F\eta + \norm{\eta}^2\xi = \lambda\,\xi\,,\quad F\eta = -\lambda\,\xi\,,
	\]
	and thus $\lambda=\frac12\,\norm{\eta}^2$. But then
	\[
		\det(F+\xi\otimes\eta) = \det F + \iprod{\Cof(F)\eta,\,\xi} = (1+\iprod{\eta,\, F^{-1}\xi})\,\det F = \Big(1-\iprod{\eta,\, \tfrac{2}{\norm{\eta}^2}\,\eta} \Big)\,\det F = -\det F\,,
	\]
	which contradicts the assumption $F,F+\xi\otimes\eta\in{\rm GL}^+(3)$.
}}
that
\begin{equation}
\widehat F\,\widehat F^T = F\,F^T\,,\qquad\widehat F = F + \xi\otimes\eta\,,\quad\det \widehat F,\; \det F > 0\qquad\implies\qquad\xi\otimes\eta = 0\,.\label{eq_one}
\end{equation}
We note that $\widehat{F}$ and $F$ are twins \cite[Sect.\ 2.5]{Silhavy_book} since they are rank-one connected and their principal stretches coincide. Here, not only their principal stretches coincide, but the left-stretch tensor is the same as well.
\begin{proof}
Since $\widehat F\,\widehat F^T = F\,F^T$, we see that $(\det \widehat F)^2 = (\det F)^2$, and by assumption \eqref{eq_one}$_3$ we can conclude that $\det \widehat F = \det F$. Since
\begin{equation}\label{eq_cross}
\det(F + \xi\otimes\eta)\ =\ \det\bigl(F(\id + F^{-1}\,\xi\otimes\eta)\bigr)\ =\ \det F \cdot \det(\id + F^{-1}\,\xi\otimes\eta)\ =\ \det F \cdot\bigl(1 + \tr(F^{-1}\,\xi\otimes\eta)\bigr)
\end{equation}
and $\det(F + \xi\otimes\eta) = \det \widehat F = \det F$, by \eqref{eq_one}$_2$ we conclude from \eqref{eq_cross}
\begin{equation}
\tr(F^{-1}\, \xi\otimes\eta)\ =\ \iprod{F^{-1}\xi\,,\,\eta}\ =\ 0\,.\label{eq_five}
\end{equation}
Assumption \eqref{eq_one}$_1$ and \eqref{eq_one}$_2$ together imply
\begin{equation}
\widehat F\,\widehat F^T\ =\ F\,F^T + F\eta\otimes\xi + \xi\otimes F\eta + \norm{\eta}^2(\xi\otimes\xi)\ =\ F\,F^T\,,\label{eq_four}
\end{equation}
thus we must have
\begin{equation}
F\eta\otimes\xi + \xi\otimes F \eta + \norm{\eta}^2(\xi\otimes\xi)= 0\,.\label{eq_six}
\end{equation}
We introduce $\widehat \xi = F^{-1}\xi$\,,\,\ $\xi = F\widehat\xi$ and insert into \eqref{eq_five} and \eqref{eq_six} to yield
\begin{equation}
F\eta \otimes F\widehat\xi + F\widehat\xi \otimes F\eta + \norm{\eta}^2\,(F\widehat\xi\otimes F\widehat\xi)\ =\ 0\,,\qquad \iprod{\widehat\xi\,,\,\eta}\ =\ 0\,.\label{eq_seven}
\end{equation}
This is equivalent to
\begin{equation}
F\,\bigl\{\,\eta\otimes\widehat\xi + \widehat\xi\otimes\eta + \norm{\eta}^2\,(\widehat\xi\otimes\widehat\xi)\,\bigr\}F^T\ =\ 0\,,\qquad\iprod{\widehat\xi\,,\,\eta}\ =\ 0\,.\label{eq_eight}
\end{equation}
Since $\det F>0$ we have as well
\begin{equation}
\eta\otimes\widehat \xi + \widehat\xi\otimes\eta + \norm{\eta}^2\,(\widehat\xi\otimes\widehat\xi)\ =\ 0\,,\qquad\iprod{\widehat\xi\,,\,\eta}\ =\ 0\,.\label{eq_nine}
\end{equation}
Multiplying \eqref{eq_nine} with $\eta\neq 0$ we obtain\ \ $\eta\underbrace{\iprod{\widehat\xi\,,\,\eta}}_{=0} + \widehat\xi\,\norm{\eta}^2 + \norm{\eta}^2\,\widehat\xi\,\underbrace{\iprod{\widehat\xi\,,\,\eta}}_{=0} = 0$. Hence, $\widehat\xi\norm{\eta}^2 = 0$ implies $\widehat\xi=0$.
\end{proof}

\end{footnotesize}

\end{document}